\documentclass[11pt]{article}

\usepackage[T1]{fontenc}
\usepackage[utf8]{inputenc}
\usepackage{amsmath,amssymb,amsthm}
\usepackage{mathtools}
\usepackage{enumitem}
\usepackage{xcolor}
\usepackage{bbm}

\newtheorem{thm}{Theorem}[section]
\newtheorem{prop}[thm]{Proposition}
\newtheorem{lem}[thm]{Lemma}
\newtheorem{cor}[thm]{Corollary}

\theoremstyle{definition}

\newtheorem{rem}[thm]{Remark}
\newtheorem{ex}[thm]{Example}

\newcommand{\ImH}{\operatorname{Im}}

\title{Decompositions of homomorphisms of abelian $p$-groups
	with $p^2$-bounded target}

\author{Stanisław Kasjan \and Justyna Kosakowska}

\date{}

\begin{document}
	
\maketitle

\begin{abstract}
	Let $p$ be a~prime. We classify homomorphisms $H:V_1\to V_2$
	between abelian $p$-groups that are direct sums of cyclic groups,
	under the assumption $p^2V_2=0$. We prove that every such homomorphism
	is isomorphic to a~direct sum of finite indecomposable representations
	and that this decomposition is unique up to isomorphism and reordering
	of the summands.
\end{abstract}

\section{Introduction}

The classification problem studied in this paper has its origins
in Birkhoff's work on subgroups of abelian groups \cite{birkhoff}.
In this context, one considers pairs $(V,U)$, where $V$ is an
abelian group and $U$ is a~subgroup of $V$, up to isomorphism.
Such a~pair corresponds to the canonical
epimorphism $V\longrightarrow V/U.$
This motivates the study of homomorphisms of abelian $p$-groups.

Let $p$ be a~fixed prime. Following \cite{kss-SH1}, we denote by $\Sigma\text{-cyc}$ the category
of abelian $p$-groups that are direct sums of cyclic groups. Put
$N_n=\mathbb Z/p^n\mathbb Z$ for $n\geq1$. Thus every object $V$ of
$\Sigma\text{-cyc}$ is of the form
$
V\simeq\bigoplus_{n\geq1}N_n^{(\kappa_n)}
$
for some cardinal numbers $\kappa_n$. We denote by
$\mathcal E_2(\Sigma\text{-cyc})$ the category whose objects are
homomorphisms
$$
H:V_1\longrightarrow V_2,
$$
where $V_1,V_2\in\Sigma\text{-cyc}$ and $p^2V_2=0$. A~morphism from
$H:V_1\to V_2$ to $H':V_1'\to V_2'$ is a~pair of homomorphisms
$(f_1,f_2)$ such that $f_1:V_1\to V_1'$, $f_2:V_2\to V_2'$ and
$$
f_2H=H'f_1.
$$
Thus $\mathcal E_2(\Sigma\text{-cyc})$ is a~subcategory of the category of
representations of the quiver $1\to2$.
We denote by $\mathcal E_2^{\mathrm{epi}}(\Sigma\text{-cyc})$
its full subcategory formed by the epimorphisms. By Kulikov's theorem,
every subgroup of a~direct sum of cyclic $p$-groups is again a~direct
sum of cyclic groups, see \cite[Theorem~8.5]{krylov-tuganbaev},
\cite[Theorem~13]{kaplansky}. In particular, $\ker H$ is a~direct sum
of cyclic groups for every object $H:V_1\to V_2$ of
$\mathcal E_2(\Sigma\text{-cyc})$.

In \cite[Theorem~1.1]{kss-SH1} it is proved that every pair $(V,U)$,
where $V$ is a~direct sum of cyclic $p$-groups and $p(V/U)=0$,
is a~direct sum of pairs $(V_i,U_i)$ with $V_i$ cyclic group. In terms of
epimorphisms, this result gives a~decomposition into maps
$N_n\to0$ and $N_n\to N_1$. We prove a~similar decomposition result
for the objects of $\mathcal E_2(\Sigma\text{-cyc})$.

The bound on the target is important. The classification of pairs $(U,V)$ under the weaker assumption, namely, $p^3U=0$, is considered infeasible, even with some restriction on $V$ ($p^8V=0$),  see
\cite[Remark~1.4(3)]{kss-SH1},
\cite{ringel-schmidmeier,schmidmeier-bounded}.
In the corresponding problem for
finite-dimensional $k[T]$-modules, the category of submodule
embeddings with $T^3$-bounded submodule and $T^8$-bounded ambient
module has wild representation type
\cite[Section~3.7]{schmidmeier-bounded}. More general Birkhoff type
problems for nilpotent linear operators, including criteria for
tame and wild representation type, are studied in
\cite{simson-birkhoff}.

The main result of the paper is the following fact.

\begin{thm}\label{thm:main}
	Every object $H:V_1\to V_2$ of
	$\mathcal E_2(\Sigma\text{-cyc})$ is isomorphic to a~direct sum,
	unique up to isomorphism and reordering of the summands, of
	representations of the following types:
	$$
	\begin{aligned}
		I_n^0&:\; N_n\to0
		&& (n\geq1),\\
		I_n^1&:\; N_n\to N_1,\quad 1\mapsto1
		&& (n\geq1),\\
		I_n^2&:\; N_n\to N_2,\quad 1\mapsto1
		&& (n\geq2),\\
		B_{m,n}^2&:\; N_m\oplus N_n\to N_2,\quad (x,y)\mapsto x+py
		&& (m\geq n+2,\ n\geq1),\\
		J_n&:\; N_n\to N_2,\quad 1\mapsto p
		&& (n\geq1),\\
		Z^1&:\; 0\to N_1\\
		Z^2&:\; 0\to N_2.
	\end{aligned}
	$$
\end{thm}

\begin{cor}\label{cor:epi}
	Every object $H:V_1\to V_2$ of
	$\mathcal E_2^{\mathrm{epi}}(\Sigma\text{-cyc})$ is isomorphic to
	a~direct sum, unique up to isomorphism and reordering of the
	summands, of representations of the types
	$I_n^0$, $I_n^1$, $I_n^2$ and $B_{m,n}^2$.
\end{cor}

\begin{rem}
	\begin{enumerate}
		\item To prove Theorem~\ref{thm:main}, it is enough to prove the existence
		of the stated decomposition. Indeed, each of the representations
		$I_n^0$, $I_n^1$, $I_n^2$, $B_{m,n}^2$, $J_n$, $Z^1$ and $Z^2$ has a~local endomorphism
		ring. Hence uniqueness up to isomorphism and reordering of the
		summands follows from the Krull--Remak--Schmidt--Azumaya theorem,
		see \cite{azumaya,facchini}.
		\item In \cite{kasjan-kosakowska-dual} we prove a~decomposition
		theorem in the setting, where $p^2\ker H=0$.
		This result requires additional assumptions and uses
		entirely different methods of proof.
	\end{enumerate}
\end{rem}

Let $H:V_1\to V_2$ be an object of
$\mathcal E_2(\Sigma\text{-cyc})$. We write
$$
V_1=\bigoplus_{n\geq1}T_n,\qquad
T_n\simeq N_n^{(\kappa_n)},\qquad
V_2\simeq N_1^{(\lambda_1)}\oplus N_2^{(\lambda_2)},
$$
where $\kappa_n,\lambda_1,\lambda_2$ are cardinals.
Set $M_0=E_0=0$ and  for $n\geq1$  put
$$
M_n=\bigoplus_{j=1}^nT_j,\qquad E_n=H(M_n).
$$
Then $E_1\subseteq E_2\subseteq\cdots$ and $\ImH H=\bigcup_{n\geq1}E_n$, because $V_1=\bigcup_{n\geq1}M_n$.

The paper and the proof of Theorem \ref{thm:main} is organized as follows. First, in Section \ref{sec-filtration}, we construct a~cyclic decomposition $V_2=\bigoplus_{\gamma\in\Gamma}W_\gamma$
such that $E_n=\bigoplus_{\gamma\in\Gamma}(E_n\cap W_\gamma)$
for every $n\geq1$. Then, in Section \ref{sec-decomposition-of-H}, we construct generators of cyclic summands
of $V_1$ whose images under $H$ lie in the summands
$W_\gamma$ and give the representations listed in
Theorem~\ref{thm:main}. The construction of these generators proceeds
by extending suitable bases of subspaces of $M_n/pM_n$. A~similar idea is used in \cite[Section~2.3]{kss-SH1}, where generators
of cyclic summands are also constructed by extending bases of
subspaces of $M_n/pM_n$. The constructions are different, however.
In \cite[Section~2.3]{kss-SH1}, the bases are chosen using the
subgroups $U_n=U\cap M_n$, whereas here they are chosen using the
decomposition of $V_2$ above and the filtration $E_n=H(M_n)$.

\section{Filtrations and cyclic decompositions}\label{sec-filtration}

\subsection{Bases and filtrations}

We first state  the following elementary fact.
We sketch the proof using more general results of \cite{gallup-sawin}.

\begin{lem}\label{lem:two-filtrations}
	Let $X$ be a~vector space over a~field $\mathbb F$ and let
	$$
	X_1\subseteq X_2\subseteq\cdots,\qquad
	Y_1\subseteq Y_2\subseteq\cdots
	$$
	be two ascending chains of subspaces of $X$.
	Then $X$ has a~basis $\mathcal B$ such that
	$$
	X_n=\operatorname{span}_{\mathbb F}(\mathcal B\cap X_n),
	\qquad
	Y_m=\operatorname{span}_{\mathbb F}(\mathcal B\cap Y_m)
	$$
	for all $n,m\geq1$.
\end{lem}

\begin{proof}
	Consider the representation
	\[
	\cdots\longleftarrow X/X_2\longleftarrow X/X_1
	\longleftarrow X
	\longrightarrow X/Y_1\longrightarrow X/Y_2\longrightarrow\cdots
	\]
	given by the natural quotient maps. This is a~representation of
	the quiver $\Theta$ of \cite[Example~1]{gallup-sawin}, whose
	underlying graph is $\mathbb{A}_{\infty,\infty}$ and all of whose
	arrows point away from the central vertex. In particular, $\Theta$
	is eventually outward, in terminology of \cite{gallup-sawin}. To a full connected subquiver $\alpha$ of $\Theta$ there corresponds a representation $X^{\alpha}$ of $\Theta$: $X^\alpha$ is one-dimensional at the vertices of $\alpha$ and zero
	elsewhere, and every arrow of $\alpha$ acts bijectively.
	By \cite[Theorem~1 and Proposition~7.1]{gallup-sawin},
	 the representation under consideration is a~direct sum of representations of type $X^\alpha$. As every arrow of $\Theta$ corresponds to an epimorphism, every subquiver $\alpha$, such that $X^{\alpha}$ appears in the decomposition, contains the central vertex of $\Theta$. The decomposition of the whole representation yields a decomposition of the space $X$. Choosing a nonzero vector in each of the one-dimensions summands one constructs a basis $\mathcal B$ of $X$. We leave it to the reader to check that this basis satisfies the condition in the lemma. 
\end{proof}

\subsection{A cyclic decomposition of the target}

The next proposition will be used to construct a~cyclic
decomposition of the target group $V_2$.

\begin{prop}\label{prop:simultaneous-target}
	Let $V\simeq N_1^{(\lambda_1)}\oplus N_2^{(\lambda_2)}$, where $\lambda_1,\lambda_2$ are cardinal numbers
	and let $E_1\subseteq E_2\subseteq\cdots\subseteq V$
	be an ascending chain of subgroups.
	Then there exists a cyclic decomposition
	\[
	V=\bigoplus_{\gamma\in\Gamma}W_\gamma,
	\]
	where every $W_\gamma$ is isomorphic to $N_1$ or $N_2$ and
	$$
	E_n=\bigoplus_{\gamma\in\Gamma}(E_n\cap W_\gamma)
	$$
	for every $n\geq1$.
	Moreover, with the convention that $n<\infty$ for every
	positive integer $n$,
\begin{enumerate}
	\item
	if $W_\gamma\simeq N_2$, then there exist  $r_\gamma\leq m_\gamma$
	in $\{1,2,\dots\}\cup\{\infty\}$
	such that
	\[
	E_n\cap W_\gamma=
	\begin{cases}
		0,&n<r_\gamma,\\
		pW_\gamma,&r_\gamma\leq n<m_\gamma,\\
		W_\gamma,&m_\gamma\leq n,
	\end{cases}
	\]
	
\item	if $W_\gamma\simeq N_1$, then there exists 
	$m_\gamma$ in $\{1,2,\dots\}\cup\{\infty\}$ such that
	$$
	E_n\cap W_\gamma=
	\begin{cases}
		0,&n<m_\gamma,\\
		W_\gamma,&m_\gamma\leq n.
	\end{cases}
	$$
\end{enumerate}	
	
\end{prop}

\begin{proof}
Put $P=pV$, $Q=V/pV$ and define
$$ \pi:V\to Q\qquad \pi(v)=v+pV,$$
	$$
	\delta:Q\longrightarrow P,
	\qquad
	\delta(v+pV)=pv.
	$$
Since $p^2V=0$, $\delta$ is well defined.
For $n\geq1$, set
$$
	B_n=\pi(E_n),
	\qquad
	A_n=E_n\cap pV,
	\qquad
	C_n=\delta^{-1}(A_n).
$$
Moreover, we put $ B_0=A_0=0, C_0=K:=\ker\delta.$ 

Applying Lemma~\ref{lem:two-filtrations} to the chains
$B_0\subseteq B_1\subseteq\cdots$ and
$C_0\subseteq C_1\subseteq\cdots$,
we obtain a~basis
$\mathcal X$ of $Q$ such that, for every $n\geq0$,	
$$B_n=\operatorname{span}(\mathcal X\cap B_n)\qquad\text{and}\qquad C_n=\operatorname{span}(\mathcal X\cap C_n).$$
For $x\in\mathcal X$,  we define
	$$
	m(x)=\min\{n\geq1:x\in B_n\}.
	$$
If in addition $x\notin K$, we define also
	$$
	r(x)=\min\{n\geq1:x\in C_n\}.
	$$
Here $\min\varnothing=\infty$.

	Since $B_n\subseteq C_n$ for every $n$, we have
	$r(x)\leq m(x)$ for every $x\in\mathcal X\setminus K$. If $m(x)<\infty$, then, since $x\in B_{m(x)}$, there exists $e_x\in E_{m(x)}$ such that
$$
x=\pi(e_x)=e_x+pV.
$$
If $m(x)=\infty$, let $e_x\in V$ be any element such that $x=\pi(e_x)$.
Fix such an element $e_x$ for every $x\in\mathcal X$. Note that if $x\notin K$, then $pe_x=\delta(x)\neq0$, and hence $e_x$ has order $p^2$. If $x\in K$, then $pe_x=0$ and, since $x\neq0$ in $Q$, the element $e_x$ has order $p$.
	
Since $K=C_0=\operatorname{span}(\mathcal X\cap K)$, the set $\{x+K:x\in\mathcal X\setminus K\}$
is a basis of $Q/K$. Moreover, since $K=\ker\delta$ and $\delta(Q)=pV$, the map $\delta$ induces an isomorphism $Q/K\xrightarrow{\sim}pV.$
Therefore
$$
\{pe_x:x\in\mathcal X\setminus K\}
$$
is a basis of $pV$.

Put $W_x=\langle e_x\rangle$. We claim that
$$
V=\bigoplus_{x\in\mathcal X}W_x.
$$
The elements $e_x+pV$, $x\in\mathcal X$, form a basis of $V/pV$, so every $v\in V$ can be written in the form
$$
v=\sum_{x\in F}a_xe_x+w,
$$
where $F\subseteq\mathcal X$ is finite and $w\in pV$. Since the elements $pe_x$, $x\in\mathcal X\setminus K$, form a basis of $pV$, we have $w\in\sum_{x\in\mathcal X}W_x$. Hence $V=\sum_{x\in\mathcal X}W_x$. \footnote{This is a case of Nakayama Lemma, see \cite[26.16]{Lam}, \cite[page 12]{Matsumura}. }
To prove that the sum is direct, suppose that
$$
\sum_{x\in F}a_xe_x=0
$$
for a finite set $F\subseteq\mathcal X$ and $a_x\in\mathbb Z$. Applying the canonical projection $\pi:V\to V/pV$, we obtain
$$
\sum_{x\in F}a_x(e_x+pV)=0.
$$
Since $V/pV$ is a vector space over $\mathbb F_p$ and the elements $e_x+pV$, $x\in\mathcal X$, are linearly independent, we have $p\mid a_x$ for every $x\in F$. For $x\in K$, this gives $a_xe_x=0$, because $e_x$ has order $p$. For $x\notin K$, write $a_x=pb_x$. Then
$$
\sum_{x\in F\setminus K}b_xpe_x=0.
$$
Since the elements $pe_x$, $x\in\mathcal X\setminus K$, are linearly independent in $pV$, each $b_x$ is divisible by $p$. As $e_x$ has order $p^2$, we again obtain $a_xe_x=0$. Thus the sum $\bigoplus_{x\in\mathcal X}W_x$ is direct.
	
We now prove that
$$
E_n=\bigoplus_{x\in\mathcal X}(E_n\cap W_x)
$$
for every $n\geq1$. Let $v\in E_n$. Then $v+pV\in B_n$. Since $B_n=\operatorname{span}(\mathcal X\cap B_n),$
we can choose integers $a_1,\ldots,a_r$ and elements
$x_1,\ldots,x_r\in\mathcal X\cap B_n$ such that
$$
v+pV=\sum_{i=1}^r a_i x_i.
$$
Since $x_i=e_{x_i}+pV$, we obtain
$v-\sum_{i=1}^r a_i e_{x_i}\in pV.$
Moreover, $x_i\in B_n$ implies $m(x_i)\leq n$, and hence $e_{x_i}\in E_{m(x_i)}\subseteq E_n.$
Since $v\in E_n$, it follows that
$v-\sum_{i=1}^r a_i e_{x_i}\in E_n.$
Therefore
$$
v-\sum_{i=1}^r a_i e_{x_i}\in E_n\cap pV=A_n.
$$

We have $C_n=\delta^{-1}(A_n)$,
$
C_n=\operatorname{span}(\mathcal X\cap C_n)$ and
$K=\operatorname{span}(\mathcal X\cap K).$
Since $\ker\delta=K$, it follows that
$$
A_n=\operatorname{span}_{\mathbb F_p}
\{pe_x:x\in\mathcal X\cap(C_n\setminus K)\}.
$$
Therefore there exist $y_1,\ldots,y_s\in
\mathcal X\cap(C_n\setminus K)$ and $b_1,\ldots,b_s\in\mathbb F_p$
such that
$
v-\sum_{i=1}^r a_i e_{x_i}
=
\sum_{j=1}^s b_jpe_{y_j}.
$
Hence
$
v=
\sum_{i=1}^r a_i e_{x_i}
+
\sum_{j=1}^s b_jpe_{y_j}.
$ For every $i$, we have $e_{x_i}\in E_n\cap W_{x_i}$,
and for every $j$, we have
$pe_{y_j}\in A_n\cap W_{y_j}\subseteq E_n\cap W_{y_j}$.
Thus $v\in\sum_{x\in\mathcal X}(E_n\cap W_x)$
and therefore
$$
E_n=\sum_{x\in\mathcal X}(E_n\cap W_x).
$$
Since $V=\bigoplus_{x\in\mathcal X}W_x$, this sum is direct.

Finally, let $x\in\mathcal X$. If $x\notin K$, then $W_x\simeq N_2$.
Note that $W_x\subseteq E_n$ if and only if $x\in B_n$
(if $x\in B_n$, then $m(x)\leq n$ and $e_x\in E_{m(x)}\subseteq E_n$),
and $pW_x\subseteq E_n$ if and only if $pe_x=\delta(x)\in A_n$,
that is, $x\in C_n$.
By the definitions of $r(x)$ and $m(x)$ we have
$$
E_n\cap W_x=
\begin{cases}
	0,&n<r(x),\\
	pW_x,&r(x)\leq n<m(x),\\
	W_x,&m(x)\leq n.
\end{cases}
$$
If $x\in K$, then $W_x\simeq N_1$ and
$$
E_n\cap W_x=
\begin{cases}
	0,&n<m(x),\\
	W_x,&m(x)\leq n.
\end{cases}
$$
This finishes the proof.
\end{proof}

\begin{ex}
	We illustrate the proof of Proposition~\ref{prop:simultaneous-target}.
	Let $V=\langle e\rangle\oplus\langle f\rangle\simeq N_2\oplus N_2$
	and put
	$$
	E_1=\langle p(e+f)\rangle,\qquad
	E_2=E_3=\langle pe,pf\rangle=pV,\qquad
	E_4=\langle e,pf\rangle,
	$$
	and $E_n=V$ for $n\geq5$.
	The original decomposition does not satisfy the conclusion
	of the proposition, since $E_1\neq0$, whereas
	$E_1\cap\langle e\rangle=E_1\cap\langle f\rangle=0$.
	
	Put $P=pV$, $Q=V/pV$, and write $q_e=e+pV$ and $q_f=f+pV$.
	The map $\delta:Q\to P$ is given by $\delta(v+pV)=pv$.
	Here $K=\ker\delta=0$. The subspaces
	$A_n=E_n\cap P$, $B_n=(E_n+P)/P$, and
	$C_n=\delta^{-1}(A_n)$ are
	$$
	\begin{aligned}
		&A_1=\langle p(e+f)\rangle,
		&&B_1=0,
		&&C_1=\langle q_e+q_f\rangle,\\
		&A_2=A_3=P,
		&&B_2=B_3=0,
		&&C_2=C_3=Q,\\
		&A_4=P,
		&&B_4=\langle q_e\rangle,
		&&C_4=Q,\\
		&A_n=P,
		&&B_n=Q,
		&&C_n=Q \qquad(n\geq5).
	\end{aligned}
	$$
	
	Set $x_1=q_e+q_f$ and $x_2=q_e$. Then
	$\mathcal X=\{x_1,x_2\}$ is a~basis of $Q$ such that
	$$
	B_n=\operatorname{span}_{\mathbb F_p}(\mathcal X\cap B_n),
	\qquad
	C_n=\operatorname{span}_{\mathbb F_p}(\mathcal X\cap C_n)
	$$
	for every $n\geq1$. Moreover, $r(x_1)=1$ and $m(x_1)=5$,
	whereas $r(x_2)=2$ and $m(x_2)=4$.
	
	We choose $e_{x_1}=e+f$ and $e_{x_2}=e$.
	Putting $W_1=\langle e+f\rangle$ and $W_2=\langle e\rangle$,
	we obtain $V=W_1\oplus W_2$ and
	$$
	\begin{aligned}
		&E_1=pW_1,\\
		&E_2=E_3=pW_1\oplus pW_2,\\
		&E_4=pW_1\oplus W_2,\\
		&E_n=W_1\oplus W_2 \qquad(n\geq5).
	\end{aligned}
	$$
\end{ex}

\section{Decomposition of the homomorphism $H:V_1\to V_2$}\label{sec-decomposition-of-H}

Applying Proposition~\ref{prop:simultaneous-target} to the chain
$E_n=H(M_n)$, we obtain a~cyclic decomposition
$$
V_2=\bigoplus_{\gamma\in\Gamma}W_\gamma,
\qquad W_\gamma\simeq N_1\ \text{or}\ N_2,
$$
such that $E_n=\bigoplus_{\gamma\in\Gamma}(E_n\cap W_\gamma)$
for every $n\geq1$.
For $W_\gamma\simeq N_2$, let $r_\gamma\leq m_\gamma$
be the elements of $\{1,2,\dots\}\cup\{\infty\}$ such that
$$
E_n\cap W_\gamma=
\begin{cases}
	0,&n<r_\gamma,\\
	pW_\gamma,&r_\gamma\leq n<m_\gamma,\\
	W_\gamma,&m_\gamma\leq n.
\end{cases}
$$
For $W_\gamma\simeq N_1$, let $m_\gamma$ be the element of $\{1,2,\dots\}\cup\{\infty\}$
such that
$$
E_n\cap W_\gamma=
\begin{cases}
	0,&n<m_\gamma,\\
	W_\gamma,&m_\gamma\leq n.
\end{cases}
$$
For every $\gamma\in\Gamma$, fix a~generator $w_\gamma$
of $W_\gamma$. Put
$$
\begin{aligned}
	&\Gamma^{\mathrm{soc}}
	=\{\gamma\in\Gamma:
	W_\gamma\simeq N_2,\ r_\gamma<m_\gamma\},\\
	& \Gamma^{\mathrm{top}}
	=\{\gamma\in\Gamma:
	W_\gamma\simeq N_2,\ m_\gamma<\infty\},\\
	&\Gamma^{(1)}
	=\{\gamma\in\Gamma:
	W_\gamma\simeq N_1\,\text{and}\,m_{\gamma}<\infty\}.
\end{aligned}
$$
Note that $\Gamma^{\mathrm{soc}}$ contains also the indices
	$\gamma$ with $m_\gamma=\infty$. We are going to construct a~cyclic decomposition  
$$
V_1=\bigoplus_{\gamma\in\Gamma^{soc}}\langle a_{\gamma}\rangle\oplus \bigoplus_{\gamma\in\Gamma^{top}}\langle b_{\gamma}\rangle \oplus \bigoplus_{\gamma\in\Gamma^{(1)}}\langle c_{\gamma}\rangle\oplus \bigoplus_{\xi\in\Xi}\langle k_{\xi}\rangle,
$$
where $\Xi$ is a set, which satisfies the following conditions.

\begin{enumerate}
\item[(i)] $\langle a_{\gamma}\rangle \simeq N_{r_{\gamma}}$ and $H(a_{\gamma})=pw_{\gamma}$ for $\gamma\in \Gamma^{\mathrm{soc}}$. 
\item[(ii)] $\langle b_{\gamma}\rangle \simeq N_{m_{\gamma}}$ and $H(b_{\gamma})=w_{\gamma}$ for $\gamma\in \Gamma^{\mathrm{top}}$. 
\item[(iii)] $\langle c_{\gamma}\rangle \simeq N_{m_{\gamma}}$ and $H(c_{\gamma})=w_{\gamma}$ for $\gamma\in \Gamma^{\mathrm{(1)}}$. 
\item[(iv)] $H(k_\xi)=0$ for $\xi\in\Xi$.
\end{enumerate}
Note that the sets $\Gamma^{\mathrm{soc}}$ and
	$\Gamma^{\mathrm{top}}$ need not be disjoint: for $r_\gamma<m_\gamma<\infty$ both
	$a_\gamma$ and $b_\gamma$ occur.
	
\subsection{Construction of the generators}

We construct these generators by induction on $n$, extending
at each step the cyclic decomposition of $M_{n-1}$ to a~cyclic decomposition
of $M_n$.
Given a natural number $n$ we set
$$
\begin{aligned}
	&\Gamma_n^{\mathrm{soc}}
	=\{\gamma\in\Gamma:
	W_\gamma\simeq N_2,\ r_\gamma=n<m_\gamma\},\\
	&\Gamma_n^{\mathrm{top}}
	=\{\gamma\in\Gamma:
	W_\gamma\simeq N_2,\ m_\gamma=n\},\\
	&\Gamma_n^{(1)}
	=\{\gamma\in\Gamma:
	W_\gamma\simeq N_1,\ m_\gamma=n\}.
\end{aligned}
$$
Note that $\Gamma_n^{\mathrm{soc}}$ contains also the indices
$\gamma$ with $r_\gamma=n$ and $m_\gamma=\infty$. Note also that $\Gamma^X$ is the union of all $\Gamma^X_n$ over natural $n$, where $X$ stands for one of the symbols $\mathrm{soc}$, $\mathrm{top}$, $(1)$.
	We also put $\Gamma^X_{<n}=\bigcup_{k<n}\Gamma^X_k$ and
	$\Xi_{<n}=\bigcup_{k<n}\Xi_k$, where $\Xi_k$ is the index set of the
	generators $k_\xi$ chosen in the $k$-th step below. 
We start with $M_0=0$ and the empty set of generators.

Let $n\geq1$ and suppose that we constructed a~cyclic
decomposition

$$
		M_{n-1}=\bigoplus_{\gamma\in\Gamma_{<n}^{\mathrm{soc}}}\langle a_{\gamma}\rangle\oplus \bigoplus_{\gamma\in\Gamma_{<n}^{\mathrm{top}}}\langle b_{\gamma}\rangle \oplus \bigoplus_{\gamma\in\Gamma_{<n}^{(1)}}\langle c_{\gamma}\rangle\oplus \bigoplus_{\xi\in\Xi_{<n}}\langle k_{\xi}\rangle,
		$$
		and the conditions (i)-(iv), formulated above, are satisfied (with $\Gamma^X$ and $\Xi$ replaced by $\Gamma^X_{<n}$ and $\Xi_{<n}$, respectively).
	Let us denote by $\mathcal G_{n-1}$ a set which is in bijection with the disjoint union of $\Gamma_{<n}^{\mathrm{soc}}$, $\Gamma_{<n}^{\mathrm{top}}$, $\Gamma_{<n}^{(1)}$ and $\Xi_{<n}$, then the above decomposition can be shortly written as $$
		M_{n-1}=\bigoplus_{g\in\mathcal G_{n-1}}\langle g\rangle,
		$$
		where $g=a_{\gamma}$ if $g$ corresponds to the element $\gamma$ in $\Gamma_{<n}^{\mathrm{soc}}$ etc.

Put
$$
A_{n-1}=\frac{M_{n-1}+pM_n}{pM_n}
\subseteq\frac{M_n}{pM_n}.
$$
Since $M_n=M_{n-1}\oplus T_n$, we have
$M_{n-1}\cap pM_n=pM_{n-1}$. By isomorphism theorem, we have
$M_{n-1}/pM_{n-1}\simeq A_{n-1}$ and
$$
\mathcal A_{n-1}=\{g+pM_n:g\in\mathcal G_{n-1}\}
$$
is a~basis of $A_{n-1}$.

We extend $\mathcal A_{n-1}$ to a~basis of $M_n/pM_n$
in two steps. Consider the linear map
$$
\varphi_n:\frac{M_n}{pM_n}\longrightarrow
\frac{E_n}{E_{n-1}+pE_n},
\qquad
\varphi_n(x+pM_n)=H(x)+(E_{n-1}+pE_n).
$$
It is well defined because $H(pM_n)=pE_n$, and it is
surjective because $H(M_n)=E_n$. Put $D_n=\ker\varphi_n$.
Then
$$
D_n=
\left\{
x+pM_n\in\frac{M_n}{pM_n}:
H(x)\in E_{n-1}+pE_n
\right\}.
$$
Since $H(M_{n-1})=E_{n-1}$, the map $\varphi_n$ vanishes
on $A_{n-1}$, so $A_{n-1}\subseteq D_n$.

\medskip
\noindent\textbf{First step of extension.}
Extend the basis $\mathcal A_{n-1}$ of $A_{n-1}$ to a basis
\[
\mathcal A_{n-1}\mathbin{\dot\cup}
\{x_\xi+pM_n:\xi\in \Xi_{n}\}
\]
of $D_n$. Since $x_\xi+pM_n\in D_n$, choose
$y_\xi\in M_{n-1}$ and $z_\xi\in M_n$ such that
\[
H(x_\xi)=H(y_\xi)+pH(z_\xi),
\]
and set $k_\xi=x_\xi-y_\xi-pz_\xi$. Then $H(k_\xi)=0$.
Moreover,
\[
(k_\xi+pM_n)-(x_\xi+pM_n)
=-y_\xi+pM_n\in A_{n-1},
\]
so the classes $k_\xi+pM_n$ and $x_\xi+pM_n$ have the same
image in $D_n/A_{n-1}$. Consequently,
\[
\mathcal A_{n-1}\mathbin{\dot\cup}
\{k_\xi+pM_n:\xi\in \Xi_n\}
\]
is a basis of $D_n$.

\medskip
\noindent\textbf{Second step of extension.}
Since $\varphi_n$ is surjective and $D_n=\ker\varphi_n$,
it remains to choose classes whose images under $\varphi_n$
form a~basis of $E_n/(E_{n-1}+pE_n)$.

\begin{lem}\label{lem:quotient-En}
	The family
	$$
	\begin{aligned}
		&\bigl\{
		pw_\gamma+(E_{n-1}+pE_n):
		\gamma\in\Gamma_n^{\mathrm{soc}}
		\bigr\}
		&\mathbin{\dot\cup}\;\;
		\bigl\{
		w_\gamma+(E_{n-1}+pE_n):
		\gamma\in\Gamma_n^{\mathrm{top}}\cup\Gamma_n^{(1)}
		\bigr\}
	\end{aligned}
	$$
	is a~basis of the $\mathbb F_p$-vector space
	$E_n/(E_{n-1}+pE_n)$.
\end{lem}

\begin{proof}
	Since $E_k=\bigoplus_{\gamma\in\Gamma}(E_k\cap W_\gamma)$
	for $k=n-1,n$, we have
	$$
	\frac{E_n}{E_{n-1}+pE_n}
	\simeq
	\bigoplus_{\gamma\in\Gamma}
	\frac{E_n\cap W_\gamma}
	{(E_{n-1}\cap W_\gamma)+p(E_n\cap W_\gamma)}.
	$$
	By Proposition~\ref{prop:simultaneous-target}, the summand
	indexed by $\gamma$ is naturally isomorphic to
	$$
	\begin{cases}
		pW_\gamma,&\gamma\in\Gamma_n^{\mathrm{soc}},\\
		W_\gamma/pW_\gamma,&\gamma\in\Gamma_n^{\mathrm{top}},\\
		W_\gamma,&\gamma\in\Gamma_n^{(1)},\\
		0,&\text{otherwise}.
	\end{cases}
	$$
	Each nonzero summand is one-dimensional over $\mathbb F_p$,
	generated by $pw_\gamma$, $w_\gamma+pW_\gamma$, or $w_\gamma$,
	respectively. These generators correspond to the classes
	listed in the statement, which therefore form a~basis.
\end{proof}

For $\gamma\in\Gamma_n^{\mathrm{soc}}$, we have
$pw_\gamma\in E_n=H(M_n)$, so we can choose
$a_\gamma\in M_n$ with $H(a_\gamma)=pw_\gamma$.
For $\gamma\in\Gamma_n^{\mathrm{top}}\cup\Gamma_n^{(1)}$,
we have $w_\gamma\in E_n$. We therefore choose
$b_\gamma\in M_n$ with $H(b_\gamma)=w_\gamma$
for $\gamma\in\Gamma_n^{\mathrm{top}}$, and
$c_\gamma\in M_n$ with $H(c_\gamma)=w_\gamma$
for $\gamma\in\Gamma_n^{(1)}$.

By Lemma~\ref{lem:quotient-En}, the images under $\varphi_n$
of these new classes form a~basis of $E_n/(E_{n-1}+pE_n)$.
Together with the basis of $D_n$ constructed in the first step,
they give the following basis of $M_n/pM_n$:
$$
\begin{aligned}
	\mathcal B_n={}&
	\mathcal A_{n-1}
	\mathbin{\dot\cup}
	\{k_\xi+pM_n:\xi\in \Xi_n\}\\
	&\mathbin{\dot\cup}
	\{a_\gamma+pM_n:\gamma\in\Gamma_n^{\mathrm{soc}}\}\\
	&\mathbin{\dot\cup}
	\{b_\gamma+pM_n:\gamma\in\Gamma_n^{\mathrm{top}}\}\\
	&\mathbin{\dot\cup}
	\{c_\gamma+pM_n:\gamma\in\Gamma_n^{(1)}\}.
\end{aligned}
$$
The generators already chosen in $M_{n-1}$ remain unchanged.

The following consequence of \cite[Lemma~2.2 and Proposition 2.3]{kss-SH1}
allows us to pass from a~basis of $M_n/pM_n$
to a~cyclic decomposition of $M_n$.

\begin{prop}\label{prop-lifting} 	Let $M_n=M_{n-1}\oplus T_n$, where $p^{n-1}M_{n-1}=0$
	and $T_n$ is a~direct sum of copies of $N_n$. Suppose that
	$$
	M_{n-1}=\bigoplus_{g\in\mathcal G_{n-1}}\langle g\rangle.
	$$
		Let $\mathcal C_n$ be a set of elements of $M_n$.
	Suppose that the classes $g+pM_n$, indexed by
	$g\in\mathcal G_{n-1}$, together with the classes
	$c+pM_n$, $c\in\mathcal C_n$, form a~basis of $M_n/pM_n$.
	Then every element of $\mathcal C_n$ has order $p^n$ and
	$$
	M_n=
	\bigoplus_{g\in\mathcal G_{n-1}}\langle g\rangle
	\oplus
	\bigoplus_{c\in \mathcal C_n}\langle c\rangle=M_{n-1}\oplus
	\bigoplus_{c\in \mathcal C_n}\langle c\rangle.
	$$
\end{prop}

Applying Proposition~\ref{prop-lifting} to the set
$\mathcal C_n$ of the new generators $k_\xi$ ($\xi\in\Xi_n$), $a_\gamma$
	($\gamma\in\Gamma^{\mathrm{soc}}_n$), $b_\gamma$ ($\gamma\in\Gamma^{\mathrm{top}}_n$) and
	$c_\gamma$ ($\gamma\in\Gamma^{(1)}_n$),
we obtain a~cyclic decomposition of $M_n$ extending that
of $M_{n-1}$. Each new generator has order $p^n$
and the required image under $H$. This completes
the induction.

\subsection{Completion of the proof of Theorem \ref{thm:main}}

By Proposition~\ref{prop-lifting}, each step yields a cyclic
decomposition of $M_n$ extending that of $M_{n-1}$, with all new
generators having order $p^n$. Since $V_1=\bigcup_{n\geq1}M_n$ and the generators chosen
at each step remain unchanged at all later steps, we obtain
a~cyclic decomposition of $V_1$.

Together with $V_2=\bigoplus_{\gamma\in\Gamma}W_\gamma$ and the prescribed images of the
generators, this decomposition of $V_1$, with $\Xi=\bigcup_{n\geq1}\Xi_n$, expresses $H$ as
a~direct sum of the maps $\langle k_\xi\rangle\to0$, $\xi\in\Xi$, and of the maps
$$
H_\gamma:V_1^\gamma\longrightarrow W_\gamma,\qquad \gamma\in\Gamma,
$$
where $V_1^\gamma$ is generated by the generators $a_\gamma,b_\gamma,c_\gamma$ attached to
$\gamma$ (if there are any). For $\xi\in\Xi_n$ the map $\langle k_\xi\rangle\to0$ is $I_n^0$.
The set $\Gamma$ is the disjoint union of the following subsets, and for each of them we
identify $H_\gamma$.
\begin{enumerate}
	\item[(a)] $\gamma\in\Gamma^{(1)}$: $H_\gamma:\langle c_\gamma\rangle\simeq N_{m_\gamma}\to
	W_\gamma\simeq N_1$, $c_\gamma\mapsto w_\gamma$, that is, $I^1_{m_\gamma}$.
	\item[(b)] $\gamma\in\Gamma^{\mathrm{top}}\setminus\Gamma^{\mathrm{soc}}$, i.e.\
	$r_\gamma=m_\gamma<\infty$: $H_\gamma:\langle b_\gamma\rangle\simeq N_{m_\gamma}\to
	W_\gamma\simeq N_2$, $b_\gamma\mapsto w_\gamma$, that is, $I^2_{m_\gamma}$.
	\item[(c)] $\gamma\in\Gamma^{\mathrm{soc}}\setminus\Gamma^{\mathrm{top}}$, i.e.\
	$r_\gamma<m_\gamma=\infty$: $H_\gamma:\langle a_\gamma\rangle\simeq N_{r_\gamma}\to
	W_\gamma\simeq N_2$, $a_\gamma\mapsto pw_\gamma$, that is, $J_{r_\gamma}$.
	\item[(d)] $\gamma\in\Gamma^{\mathrm{soc}}\cap\Gamma^{\mathrm{top}}$, i.e.\
	$r_\gamma<m_\gamma<\infty$:
	$$
	H_\gamma:\langle b_\gamma\rangle\oplus\langle a_\gamma\rangle\simeq
	N_{m_\gamma}\oplus N_{r_\gamma}\longrightarrow W_\gamma\simeq N_2,
	\qquad b_\gamma\mapsto w_\gamma,\quad a_\gamma\mapsto pw_\gamma .
	$$
	If $m_\gamma\geq r_\gamma+2$, this is the bipicket $B^2_{m_\gamma,r_\gamma}$.
	If $m_\gamma=r_\gamma+1$, put $a_\gamma'=a_\gamma-pb_\gamma$. Then $H(a_\gamma')=0$,
	$a_\gamma'$ has order $p^{r_\gamma}$ and
	$\langle a_\gamma\rangle\oplus\langle b_\gamma\rangle
	=\langle a_\gamma'\rangle\oplus\langle b_\gamma\rangle$,
	so $H_\gamma\simeq I^0_{r_\gamma}\oplus I^2_{m_\gamma}$.
	\item[(e)] $W_\gamma\simeq N_2$ and $r_\gamma=\infty$, or $W_\gamma\simeq N_1$ and
	$m_\gamma=\infty$: no generator is attached to $\gamma$, so $V_1^\gamma=0$ and $H_\gamma$
	is $Z^2$, respectively $Z^1$. Indeed, in this case
	$W_\gamma\cap\ImH H=\bigcup_{n\geq1}(E_n\cap W_\gamma)=0$.
\end{enumerate}
This gives the required decomposition of $H$.

\section{Examples}

\begin{ex}
	Let
	\[
	V_1=\bigoplus_{i\geq1}
	\bigl(\langle a_i\rangle\oplus\langle b_i\rangle\bigr),
	\qquad
	V_2=\bigoplus_{i\geq1}\langle e_i\rangle,
	\]
	where $a_i$, $b_i$, and $e_i$ have orders $p^i$, $p^{i+2}$,
	and $p^2$, respectively. Set $e_0=0$ and define
	\[
	H(a_i)=pe_i,\qquad H(b_i)=e_i-e_{i-1}.
	\]
	Then $H$ is surjective, since $H(b_1+\cdots+b_i)=e_i$.
	We follow the construction in the proof of
	Theorem~\ref{thm:main}. All sums over empty index sets are zero.
	
	\medskip
	\noindent\textbf{The target decomposition.}
	We have
	\[
	M_n=
	\bigoplus_{1\leq i\leq n}\langle a_i\rangle
	\oplus
	\bigoplus_{1\leq i\leq n-2}\langle b_i\rangle,
	\qquad
	E_n=
	\bigoplus_{1\leq i\leq n-2}\langle e_i\rangle
	\oplus
	\bigoplus_{\substack{i\geq1\\n-2<i\leq n}}p\langle e_i\rangle.
	\]
		Thus the decomposition $V_2=\bigoplus_{i\geq1}\langle e_i\rangle$ already satisfies the
	conclusion of Proposition~\ref{prop:simultaneous-target}, and we may take
	$W_i=\langle e_i\rangle$, with $r_i=i$ and $m_i=i+2$. Consequently,
	\[
	\Gamma_n^{\mathrm{soc}}=\{n\},\qquad
	\Gamma_n^{\mathrm{top}}=
	\begin{cases}
		\varnothing,&n=1,2,\\
		\{n-2\},&n\geq3,
	\end{cases}
	\qquad
	\Gamma_n^{(1)}=\varnothing.
	\]
	
	\medskip
	\noindent\textbf{The first step of extension.}
	For $n\geq3$, we have
	$M_n=M_{n-1}\oplus\langle a_n\rangle
	\oplus\langle b_{n-2}\rangle$.
	Hence the classes of $a_n$ and $b_{n-2}$ form a~basis of
	$M_n/(M_{n-1}+pM_n)$. Their images under the map induced by
	$\varphi_n$ are the classes of $pe_n$ and $e_{n-2}$,
	which form a~basis of $E_n/(E_{n-1}+pE_n)$.
	For $n=1,2$, only the class of $a_n$ occurs.
	Therefore $D_n=A_{n-1}$, and no new generators $k_\xi$
	are needed.
	
	\medskip
	\noindent\textbf{The second step of extension.}
	For the index $n\in\Gamma_n^{\mathrm{soc}}$, we choose $a_n$,
	since $H(a_n)=pe_n$.
	Set $\beta_0=0$ and
	\[
	\beta_i=b_i+\beta_{i-1}=b_1+\cdots+b_i.
	\]
For $n\geq3$, the index $n-2\in\Gamma_n^{\mathrm{top}}$
		requires a~preimage of $w_{n-2}=e_{n-2}$. The element $b_{n-2}$ does not work, since
		$H(b_{n-2})=e_{n-2}-e_{n-3}$, so we choose $\beta_{n-2}$, for which
		$H(\beta_{n-2})=e_{n-2}$. Proposition~\ref{prop-lifting} gives
	\[
	M_n=
	\bigoplus_{1\leq i\leq n}\langle a_i\rangle
	\oplus
	\bigoplus_{1\leq i\leq n-2}\langle\beta_i\rangle.
	\]
	
	\medskip
	\noindent\textbf{The resulting decomposition.}
	Taking the union over $n$, we obtain
	\[
	V_1=\bigoplus_{i\geq1}
	\bigl(\langle\beta_i\rangle\oplus\langle a_i\rangle\bigr),
	\qquad
	H(\beta_i)=e_i,\quad H(a_i)=pe_i.
	\]
	Since $\beta_i$ has order $p^{i+2}$ and $a_i$ has order $p^i$,
	each pair gives the bipicket $B_{i+2,i}^2$. Hence
	\[
	H\simeq\bigoplus_{i\geq1}B_{i+2,i}^2.
	\]
\end{ex}

\begin{ex}
	We illustrate the case where $H$ is not surjective.
	Let $V_1=\langle a\rangle\simeq N_1$,
	$V_2=\langle e\rangle\oplus\langle f\rangle\simeq N_2\oplus N_1$,
	and define $H(a)=pe+f$. Then $M_n=V_1$ and
	$E_n=\langle pe+f\rangle$ for all $n\geq1$. The decomposition
	$V_2=\langle e\rangle\oplus\langle f\rangle$ does not satisfy the
	conclusion of Proposition~\ref{prop:simultaneous-target}, since
	$E_1\cap\langle e\rangle=E_1\cap\langle f\rangle=0$.

	Write $q_e=e+pV_2$ and $q_f=f+pV_2$. Then $C_0=K=\langle q_f\rangle$,
	and for $n\geq1$ we have
	$B_n=\langle q_f\rangle$, $A_n=E_n\cap pV_2=0$ and $C_n=K$.
	The basis $\mathcal X=\{q_f,q_e\}$ of $Q$ satisfies the conditions
	of Lemma~\ref{lem:two-filtrations}, with $q_f\in K$, $m(q_f)=1$,
	and $q_e\notin K$, $r(q_e)=m(q_e)=\infty$.
		Since $m(q_f)=1$, the preimage $e_{q_f}$ of $q_f$ must lie in $E_1$, so we take
	$e_{q_f}=pe+f$. 	Since $m(q_e)=\infty$, the proof of Proposition~\ref{prop:simultaneous-target} imposes no
	condition on the preimage of $q_e$, and we take $e_{q_e}=e$.
		Thus the new decomposition differs from $V_2=\langle e\rangle\oplus\langle f\rangle$ only
	in that $f$ is replaced by $pe+f$.
	Hence
	$$
	V_2=W_1\oplus W_2,\qquad
	W_1=\langle pe+f\rangle\simeq N_1,\quad
	W_2=\langle e\rangle\simeq N_2,
	$$
	and $E_n=W_1$ for all $n\geq1$. Here $\Gamma_1^{(1)}=\{1\}$, and
	we choose $c_1=a$. Consequently
	$$
H\simeq I_1^1\oplus Z^2,
$$
where $Z^2$ comes from $W_2$, for which $r_2=m_2=\infty$.\end{ex}

{\bf Acknowledgements.} ChatGPT (OpenAI, GPT-5.6 Sol), accessed through a ChatGPT Plus subscription, was used to assist with language editing and to review the presentation, consistency, and correctness of some mathematical arguments. All suggestions were critically assessed by the authors, who take full responsibility for the mathematical content of the paper.

\end{document}